\documentclass[11pt,letterpaper]{article}
\usepackage{graphicx}
\usepackage{epsfig,subfigure}
\usepackage{amssymb}
\usepackage{amsthm}
\usepackage{amsfonts}
\usepackage{amsmath}
\usepackage[margin=1in]{geometry}

\newtheorem{theorem}{Theorem}
\newtheorem{lemma}{Lemma}
\newtheorem{corollary}[theorem]{Corollary}
\newtheorem{proposition}[theorem]{Proposition}
\newtheorem{definition}{Definition}

\newtheorem{example}[theorem]{Example}

\newcommand{\td}{\mathbb{TD}}
\newcommand{\R}{{\mathbb R}}

\newcommand{\bx}{\mathbf{x}}
\newcommand{\ba}{\mathbf{a}}
\newcommand{\bb}{\mathbf{b}}
\newcommand{\bh}{\mathbf{h}}
\newcommand{\bX}{\mathbf{X}}
\newcommand{\bxa}{\mathbf{x}^\ast}
\newcommand{\by}{\mathbf{y}}
\newcommand{\bz}{\mathbf{z}}
\newcommand{\bs}{\mathcal{S}}
\newcommand{\bc}{\mathcal{C}}
\newcommand{\bk}{\mathcal{K}}
\newcommand{\ones}{\mathbf{1}}
\newcommand{\bd}{\mathbf{d}}
\newcommand{\bdelta}{\boldsymbol\delta}
\newcommand{\bg}{\mathbf{g}}

\newcommand{\bm}{\mathbf{m}}
\newcommand{\G}{g}
\renewcommand{\S}{\mathbb{S}}
\newcommand{\tz}{\tilde{\mathbf{\bz}}}

\title{Computational and Statistical Tradeoffs via Convex Relaxation}
\author{Venkat Chandrasekaran$^c$ and Michael I. Jordan$^b$ \thanks{Email: venkatc@caltech.edu, jordan@cs.berkeley.edu} \vspace{0.25in} \\ $^c$ Department of Computing and Mathematical Sciences \\ California Institute of Technology \\ Pasadena, CA 91125 USA  \vspace{0.1in} \\ $^b$ Departments of Statistics and of Electrical Engineering and Computer Sciences \\ University of California -- Berkeley \\ Berkeley, CA 94720 USA}
\date{November 26, 2012}

\begin{document}

\maketitle

\begin{abstract}
In modern data analysis, one is frequently faced with statistical inference problems involving massive datasets.  Processing such large datasets is usually viewed as a substantial computational challenge.  However, if data are a statistician's main resource then access to more data should be viewed as an asset rather than as a burden.  In this paper we describe a computational framework based on convex relaxation to reduce the computational complexity of an inference procedure when one has access to increasingly larger datasets.  Convex relaxation techniques have been widely used in theoretical computer science as they give tractable approximation algorithms to many computationally intractable tasks.  We demonstrate the efficacy of this methodology in statistical estimation in providing concrete time-data tradeoffs in a class of denoising problems.  Thus, convex relaxation offers a principled approach to exploit the statistical gains from larger datasets to reduce the runtime of inference algorithms.
\end{abstract}

\textbf{Keywords}: massive datasets; high-dimensional statistics; convex relaxation; convex geometry

\section*{Introduction}

The rapid growth in the size and scope of datasets in science and technology has
created a need for novel foundational perspectives on data analysis that blend
computer science and statistics.  That classical perspectives from these fields
are not adequate to address emerging problems in ``Big Data'' is apparent
from their sharply divergent nature at an elementary level---in computer
science, the growth of the number of data points is a source of ``complexity''
that must be tamed via algorithms or hardware, whereas in statistics, the growth
of the number of data points is a source of ``simplicity'' in that inferences
are generally stronger and asymptotic results can be invoked.  In classical
statistics, where one considers the increase in inferential accuracy as the
number of data points grows, there is little or no consideration of computational
complexity.  Indeed, if one imposes the additional constraint---prevalent in
real-world applications---that a certain level of inferential accuracy be achieved
within a limited time budget, classical theory provides no guidance as to how to design
an inferential strategy.\footnote{Note that classical statistics contains a branch known as \emph{sequential
analysis} that does discuss methods that stop collecting data points after a target error level has been reached (see, e.g., \cite{Lai2001}), but this is different from the computational complexity guarantees (the number of steps that a computational procedure requires) that are our focus.}  In classical computer science, practical solutions to large-scale problems are often framed in terms of approximations to idealized
problems, but even when such approximations are sought, they are rarely expressed
in terms of the coin of the realm of the theory of inference---the statistical
risk function.  Thus there is little or no consideration of the idea that
computation can be simplified in large datasets because of the enhanced
inferential power in the data.  In general, in computer science, datasets
are not viewed formally as a resource on a par with time and space
(such that the more of the resource the better).

On intuitive grounds it is not implausible that strategies can be designed that
yield monotonically improving risk as data accumulate, even in the face of a
time budget.  In particular, if an algorithm simply ignores all future data once
a time budget is exhausted, then statistical risk will not increase (under various
assumptions that may not be desirable in practical applications).  Alternatively,
one might allow linear growth in the time budget (for example, in a real-time
setting), and attempt to achieve such growth via a subsampling strategy where
some fraction of the data are dropped.  Executing such a strategy may be
difficult, however, in that the appropriate fraction depends on the risk
function and thus on a mathematical analysis that may be difficult to carry out.
Moreover, subsampling is a limited strategy for controlling computational
complexity.  More generally, one would like to consider some notion of
``algorithm weakening,'' where as data accumulate one can back off to
simpler algorithmic strategies that nonetheless achieve a desired risk.
The challenge is to do this in a theoretically sound manner.

We base our approach to this problem on the notion of a ``time-data complexity class.''
In particular, we define a class $\td(t(p),n(p),\epsilon(p))$ of parameter estimation
problems in which a $p$-dimensional parameter underlying an unknown population can
be estimated with a risk of $\epsilon(p)$ given $n(p)$ i.i.d.\ samples using an inference
procedure with runtime $t(p)$.  Our definition parallels the definition of the TISP
complexity class in computational complexity theory for describing algorithmic tradeoffs
between time and space resources \cite{AroB2009}.  In this formalization, classical
results in estimation theory can be viewed as emphasizing the tradeoffs between the
second and third parameters (amount of data and risk).  Our focus in this paper is
to fix $\epsilon(p)$ to some desired level of accuracy and to investigate the tradeoffs
between the first two parameters, namely runtime and dataset size.

Although classical statistics gave little consideration to computational complexity,
computational issues have come increasingly to the fore in modern ``high-dimensional
statistics''~\cite{BuhV2011}, where the number of parameters $p$ is relatively large
and the number of data points $n$ relatively small.  In this setting, methods based
on convex optimization have been emphasized (in particular methods based on $\ell_1$
penalties).  This is due in part to the favorable analytic properties of convex functions
and convex sets, but also due to the fact that such methods tend to have favorable
computational scaling.  However, the treatment of computation has remained informal,
with no attempt to characterize tradeoffs between computation time and estimation
quality.  In our work we aim explicitly at such tradeoffs, in the setting in which
both $n$ and $p$ are large.

To develop a notion of ``algorithm weakening'' that combines computational and
statistical considerations, we consider estimation procedures for which we can
characterize the computational benefits as well as the loss in estimation performance
due to the use of weaker algorithms.  Reflecting the fact that the space of all
algorithms is poorly understood, we retain the focus on convex optimization from
high-dimensional statistics, but we consider parameterized hierarchies of
optimization procedures in which a form of algorithm weakening is obtained by
employing successively weaker outer approximations to convex sets.  Such
\emph{convex relaxations} have been widely used to give efficient approximation
algorithms for intractable problems in computer science \cite{Vaz2004}.
As we will discuss, a precise characterization of both the estimation
performance and the computational complexity of employing a particular
relaxation of a convex set can be obtained by appealing to convex geometry
and to results on the complexity of solving convex programs.  Specifically, the
tighter relaxations in these families offer better approximation quality (and in
our context better estimation performance) but such relaxations are computationally
more complex.  On the other hand the weaker relaxations are computationally more
tractable, and they can provide the same estimation performance as the tighter
ones but with access to more data.  In this manner, convex relaxations provide a
principled mechanism to weaken inference algorithms in order to reduce the runtime
in processing larger datasets.

To demonstrate explicit tradeoffs in high-dimensional, large-scale inference, we
focus for simplicity and concreteness on estimation in sequence models \cite{Joh2011}:
\begin{equation}
\by = \bxa + \sigma \bz, \label{eq:denoise}
\end{equation}
where $\sigma > 0$, the noise vector $\bz \in \R^p$ is standard normal, and the unknown parameter $\bxa$ belongs to a known subset $\bs \subset \R^p$.  The objective is to estimate $\bxa$ based on $n$ independent observations $\{\by_i\}_{i=1}^n$ of $\by$.  This denoising setup has a long history and has been at the center of some remarkable results in the high-dimensional setting over the past two decades beginning with the papers of Donoho and Johnstone \cite{Don1995,DonJ1998}.  The estimators discussed next proceed by first computing the sample mean $\bar{\by} = \sum_{i=1}^n \by_i$ and then using $\bar{\by}$ as input to a suitable convex program.  Of course, this is equivalent to a denoising problem in which the noise variance is $\sigma^2 / n$ and we are given just one sample. The reason we consider the elaborate two-step procedure is to account more accurately both for data aggregation and for subsequent processing in our runtime calculations.  Indeed, in a real-world setting one is typically faced with a massive dataset in unaggregated form, and when both $p$ and $n$ may be large, summarizing the data before any further processing can itself be an expensive computation.  As will be seen in concrete calculations of time-data tradeoffs, the number of operations corresponding to data aggregation is sometimes comparable to or even larger than the number of operations required for subsequent processing in a massive data setting.

In order to estimate $\bxa$, we consider the following natural shrinkage estimator given by a projection of the sample mean $\bar{\by}$ onto a convex set $\bc$ that is an outer approximation to $\bs$, i.e., $\bs \subset \bc$:
\begin{equation}
\hat{\bx}_n(\bc) = \arg \min_{\bx \in \R^p} ~~~ \frac{1}{2}\left\|\bar{\by} - \bx \right\|_{\ell_2}^2 ~~~ \mathrm{s.t.} ~~~ \bx \in \bc. \label{eq:shrink}
\end{equation}
We study the estimation performance of a family of shrinkage estimators $\{\hat{\bx}_n(\bc_i)\}$ that employ as the convex constraint one of a sequence of convex outer approximations $\{\bc_i\}$ with $\bc_1 \supset \bc_2 \supset \cdots \supset \bs$.  Given the same number of samples, using a weaker relaxation such as $\bc_1$ leads to an estimator with a larger risk than would result from using a tighter relaxation such as $\bc_2$.  On the other hand, given access to more data samples the weaker approximations provide the same estimation guarantees as the tighter ones.  In settings in which computing a weaker approximation is more tractable than computing a tighter one, a natural computation/sample tradeoff arises.  We characterize this tradeoff in a number of stylized examples, motivated by problems such as collaborative filtering, learning an ordering of a collection of random variables, and inference in networks.

More broadly, this paper highlights the role of computation in estimation by jointly studying both the computational and the statistical aspects of high-dimensional inference.  Such an understanding is particularly of interest in modern inferential tasks in data-rich settings.  Further, an observation from our examples on time-data tradeoffs is that in many contexts one does not need too many extra data samples in order to go from a computationally inefficient estimator based on a tight relaxation to an extremely efficient estimator based on a weaker relaxation.  Consequently, in application domains in which obtaining more data is not too expensive it may be preferable to acquire more data with the upshot being that the computational infrastructure can be relatively less sophisticated.

We should note that we investigate only one algorithm weakening mechanism, namely convex
relaxation, and one class of statistical estimation problems, namely denoising in a
high-dimensional sequence model.  There is reason to believe, however, that the principles
described in this paper are relevant more generally.  Convex-optimization-based procedures
are employed in a variety of large-scale data analysis tasks~\cite{BoyV2004,BuhV2011}, and it is
likely to be interesting to explore hierarchies of convex relaxations in such tasks.
In addition, there are a number of potentially interesting mechanisms beyond convex
relaxation for weakening inference procedures such as dimensionality reduction or other
forms of data quantization, and approaches based on clustering or coresets.  We discuss
these and other research directions in the Conclusions.

\paragraph{Related work} A number of papers have considered computational and sample
complexity tradeoffs in the setting of learning binary classifiers.  Specifically,
several authors have described settings under which speedups in running time of a
classifier learning algorithm are possible given a substantial increase in dataset
size \cite{ShaST2012,DecGR1998,Ser2000}.  In contrast, in the denoising setup considered
in this paper, several of our examples of time-data tradeoffs demonstrate significant
 computational speedups with just a constant factor increase in dataset size.
Another attempt in the binary classifier learning setting, building on earlier
work on classifier learning in data-rich problems \cite{BotB2008}, has shown that
modest improvements in runtime (of constant factors) may be possible with access to
more data by employing the stochastic gradient descent method \cite{ShaS2008}.
Time-data tradeoffs have also been characterized in Boolean network training from
time series data \cite{PerH2010}, but the computational speedups offered there are
from exponential-time algorithms to slightly faster but still exponential-time
algorithms.  Two recent papers \cite{AmiW2009,Kol2001} have considered time-data
tradeoffs in sparse principal component analysis (PCA) and in biclustering, in which
a sparse rank-one matrix is corrupted by noise and the objective is to recover
the support of the matrix.  We also study time-data tradeoffs in the estimation
of a sparse rank-one matrix, but from a denoising perspective.  In our discussion
of Example 3 below we discuss the differences between our problem setup and these
latter two papers.  As a general contrast to all these previous results, a major
contribution of the present paper is the demonstration of the efficacy of convex
relaxation as a powerful algorithm weakening mechanism for processing massive
datasets in a broad range of settings.

\paragraph{Paper outline} The main sections of this paper proceed in the following sequence.  The next section describes a framework for formally stating results on time-data tradeoffs.  Then we provide some background on convex optimization and relaxations of convex sets.  Following this we investigate in detail the denoising problem \eqref{eq:denoise}, and characterize the risk obtained when one employs a convex programming estimator of the type \eqref{eq:shrink}.  Subsequently, we give several examples of time-data tradeoffs in concrete denoising problems.  Finally, we conclude with a discussion of directions for further research.

\section*{Formally Stating Time-Data Tradeoffs}

In this section we describe a framework to state results on computational and statistical tradeoffs in estimation.  Our discussion is relevant to general parameter estimation problems and inference procedures; one may keep in mind the denoising problem \eqref{eq:denoise} for concreteness.  Consider a sequence of estimation problems indexed by the dimension $p$ of the parameter to be estimated.  Fix a risk function $\epsilon(p)$ that specifies the desired error of an estimator.  For example, in the denoising problem \eqref{eq:denoise} the error of an estimator of the form \eqref{eq:shrink} may be specified as the worst case mean squared error taken over all elements of the set $\bs$, i.e., $\sup_{\bxa \in \bs} ~ \mathbb{E}\left[\|\bxa-\hat{\bx}_n(\bc)\|^2_{\ell_2}\right]$.

One can informally view an estimation algorithm that achieves a risk of $\epsilon(p)$ by processing $n(p)$ samples with runtime $t(p)$ as a point on a two-dimensional plot such as Figure~1, with one axis representing the runtime and the other representing the sample complexity.  To be precise the axes in the plot index \emph{functions} (of $p$) that represent runtime and number of samples, but we do not emphasize such formalities and rather use these plots to provide a useful qualitative comparison of inference algorithms.  In Figure~1, procedure A requires fewer samples than procedure C to achieve the same error, but this reduction in sample complexity comes at the expense of a larger runtime.  Procedure B has both a larger sample complexity and a larger runtime than procedure C, and thus it is strictly dominated by procedure C.

\begin{figure}
\begin{center}
\epsfig{file=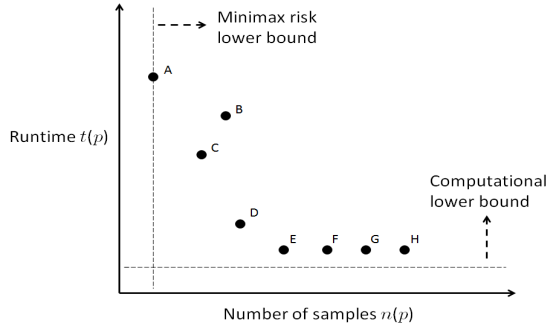,width=7.5cm,height=4.5cm} \caption{The tradeoff plot between the runtime and sample complexity in a typical parameter estimation problem.  Here the risk is assumed to be fixed to some desired level, and the points in the plot refer to different algorithms that require a certain runtime and a certain number of samples in order to achieve the desired risk.  The vertical and horizontal lines refer to lower bounds in sample complexity and in runtime, respectively.}
\end{center} \label{fig:td}
\end{figure}

Given an error function $\epsilon(p)$, there is a lower bound on the number of samples $n(p)$ required to achieve this error using any computational procedure (i.e., no constraints on $t(p)$)---such information-theoretic or minimax risk lower bounds correspond to ``vertical lines'' in the plot in Figure~1.  Characterizing these fundamental limits on sample complexity has been a traditional focus in the estimation theory literature with a fairly complete set of results available in many settings.  One can imagine asking for similar lower bounds on the computational side, corresponding to ``horizontal'' lines in the plot in Figure~1---given a desired risk $\epsilon(p)$ and access to an unbounded number of samples, what is a non-trivial lower bound on the runtime $t(p)$ of any inference algorithm that achieves a risk of $\epsilon(p)$?  Such complexity-theoretic lower bounds are significantly harder to obtain, and they remain a central open problem in computational complexity theory.

This research landscape informs the qualitative nature of the statements on time-data tradeoffs we make in this paper.  First, we will not attempt to prove combined lower bounds---as is traditionally done in the characterization of tradeoffs between physical quantities---involving $n(p)$ and $t(p)$ jointly; this is because obtaining a lower bound just on $t(p)$ remains a substantial challenge.  Hence, our time-data tradeoff results on the use of more efficient algorithms for larger datasets refer to a reduction in the \emph{upper bounds} on runtimes of estimation procedures with increases in dataset size.  Second, in any setting in which there is a computational cost associated with touching each data sample and in which the samples are exchangeable, there is a sample threshold beyond which it is computationally more efficient to throw away excess data samples than to process them in any form.  This observation suggests that there is a ``floor,'' as in Figure~1 with procedures $E,F,G$ and $H$, beyond which additional data do not lead to a reduction in runtime.  Precisely characterizing this sample threshold is in general very hard as it depends on difficult-to-obtain computational lower bounds for estimation tasks and also on the particular space of estimation algorithms that one may employ.  We will comment further on this point when we consider concrete examples of time-data tradeoffs.

In order to formally state our results concerning time-data tradeoffs, we define a resource class constrained by runtime and sample complexity as follows.
\begin{definition}
Consider a sequence of parameter estimation problems indexed by the dimension $p$ of the space of parameters that index an underlying population.  This sequence of estimation problems belongs to a \emph{time-data class $\td(t(p),n(p),\epsilon(p))$} if there exists an inference procedure for the sequence of problems with runtime upper-bounded by $t(p)$, with the number of i.i.d.\ samples processed bounded by $n(p)$, and which achieves a risk bounded by $\epsilon(p)$.
\end{definition}
We note that our definition of a time-data resource class parallels the time-space resource classes considered in complexity theory \cite{AroB2009}---in that literature $\mathrm{TISP}(t(p),s(p))$ denotes a class of problems of input size $p$ that can be solved by some algorithm using $t(p)$ operations and utilizing $s(p)$ units of space.

With this formalism, classical minimax bounds can be stated as follows.  Given some function $\bar{n}(p)$ for the number of samples, suppose a parameter estimation problem has a minimax risk of $\epsilon_{\mathrm{minimax}}(p)$ (which depends on the function $\bar{n}(p)$).  If an estimator achieving a risk of $\epsilon_{\mathrm{minimax}}(p)$ is \emph{computable} with runtime $\bar{t}(p)$, then this estimation problem lies in $\td(\bar{t}(p),\bar{n}(p),\epsilon_{\mathrm{minimax}}(p))$.  Thus the emphasis is fundamentally on the relationship between $\bar{n}(p)$ and $\epsilon_{\mathrm{minimax}}(p)$, without much focus on the computational procedure that achieves the minimax risk bound.  Our interest in this paper is to fix the risk $\epsilon(p) = \epsilon_{\mathrm{desired}}(p)$ to be equal to some desired level of accuracy, and to investigate the tradeoffs between $t(p)$ and $n(p)$ so that a parameter estimation problem lies in $\td(t(p),n(p),\epsilon_{\mathrm{desired}}(p))$.

\section*{Convex Relaxation}

In this section we describe the particular algorithmic toolbox on which we focus, namely convex programs.  Convex optimization methods offer a powerful framework for statistical inference due to the broad class of estimators that can be effectively modeled as convex programs.  Further the theory of convex analysis is useful both for characterizing the statistical properties of convex programming based estimators as well as for developing methods to compute such estimators efficiently.  Most importantly from our viewpoint, convex optimization methods provide a principled and general framework for algorithm weakening based on relaxations of convex sets.  We briefly discuss the key ideas from this literature that are relevant to this paper in this section.  A central notion to the geometric viewpoint adopted in this section is that of a \emph{convex cone}, which is a convex set that is closed under nonnegative linear combinations.

\subsection*{Representation of Convex Sets}

Convex programs refer to a class of optimization problems in which we seek to minimize a convex function over a convex constraint set \cite{BoyV2004}.  For example linear programming and semidefinite programming are two prominent subclasses in which linear functions are minimized over constraint sets given by affine spaces intersecting the nonnegative orthant (in linear programming) and the positive semidefinite cone (in semidefinite programming).  Roughly speaking convex programs are tractable to solve computationally if the convex objective function can be computed efficiently, and if membership in the convex constraint sets can be certified efficiently\footnote{More precisely, one requires an efficient \emph{separation oracle} that responds YES if the point is in the convex set, and otherwise provides a hyperplane that separates the point from the convex set.}; we will informally refer to this latter operation as computing the convex constraint set.  It is then clear that the main computational bottleneck associated with solving convex programs of the form \eqref{eq:shrink} is the efficiency of computing the constraint sets.

A central insight from the literature on convex optimization is that the complexity of computing a convex set is closely linked to how efficiently the set can be \emph{represented}.  Specifically, if a convex set can be expressed as the intersection of a small number of ``basic'' or ``elementary'' convex sets, each of which is tractable to compute, then the original convex set is also tractable to compute and one can in turn optimize over this set efficiently.  Examples of ``basic'' convex sets include affine spaces or cones such as the nonnegative orthant and the cone of positive semidefinite matrices.  Indeed, a canonical method to represent a convex set is to express the set as the intersection of a cone and an affine space.  In what follows we will consider such \emph{conic representations} of convex sets in $\R^p$.

\begin{definition}
Let $\bc \in \R^p$ be a convex set and let $\bk \in \R^p$ be a convex cone.  Then $\bc$ is said to be \emph{$\bk$-representable} if $\bc$ can be expressed as follows for $A \in \R^{m \times p}, b \in \R^m$:
\begin{equation}
\bc = \{\bx | \bx \in \bk, A \bx = b\}. \label{eq:cone-rep}
\end{equation}
Such a representation of $\bc$ is called a \emph{$\bk$-representation}.
\end{definition}

Informally, if $\bk$ is the nonnegative orthant (or the semidefinite cone) we will refer to the resulting representations as LP representations (or SDP representations), following commonly used terminology in the literature.  A virtue of conic representations of convex sets based on the orthant or the semidefinite cone is that these representations lead to a numerical recipe for solving convex optimization problems of the form \eqref{eq:shrink} via a natural associated barrier penalty \cite{NesN1995}.  The computational complexity of these procedures is polynomial in the dimension of the cone and we discuss runtimes for specific instances in our discussion of concrete examples of time-data tradeoffs.

\begin{example}
The $p$-dimensional simplex is an example of an LP representable set:
\begin{equation}
\Delta_p = \{\bx | \ones' \bx = 1, \bx \geq 0\}, \label{eq:simplex}
\end{equation}
where $\ones \in \R^p$ is the all ones vector.
\end{example}

The $p$-simplex is the set of probability vectors in $\R^p$.  The next example is one of an SDP-representable set that is commonly encountered both in optimization and in statistics.

\begin{example}
The \emph{elliptope}, or the set of correlation matrices, in the space of $m \times m$ symmetric matrices is defined as follows:
\begin{equation}
\mathcal{E}_{m \times m} = \{\bX | \bX \succeq 0, \bX_{ii} = 1 \; \forall i\}. \label{eq:elliptope}
\end{equation}
\end{example}

Conic representations are somewhat limited in their modeling capacity, and an important generalization is obtained by considering \emph{lifted} representations.  In particular the notion of \emph{lift-and-project} plays a critical role in many examples of efficient representations of convex sets.  The lift-and-project concept is simple---we wish to express a convex set $\mathcal{C} \in \R^p$ as the projection of a convex set $\mathcal{C}' \in \R^{p'}$ in some higher-dimensional space (i.e., $p' > p$).  The complexity of solving the associated convex programs is now a function of the lifting dimension $p'$.  Thus, lift-and-project techniques are useful if $p'$ is not too much larger than $p$ and if $\mathcal{C}'$ has an efficient representation in the higher-dimensional space $\R^{p'}$.  Lift-and-project provides a very powerful representation tool, as seen in the following example.

\begin{example}
The cross-polytope is the unit ball of the $\ell_1$-norm:
\begin{equation*}
B^p_{\ell_1} = \left\{ \bx \in \R^p ~ | ~ \sum_i |\bx_i| ~ \leq 1 \right\}.
\end{equation*}
The $\ell_1$-norm has been the focus of much attention recently in statistical model selection and feature selection due to its sparsity-inducing properties \cite{CanRT2006,Don2006b}.  While the cross-polytope has $2p$ vertices, a direct specification in terms of linear constraints involves $2^p$ inequalities:
\begin{equation*}
B^p_{\ell_1} = \left\{\bx \in \R^p ~|~ \sum_i \bz_i \bx_i \leq 1, ~ \forall \bz \in \{-1,+1\}^p \right\}.
\end{equation*}
However we can obtain a tractable representation by lifting to $\R^{2p}$ and then projecting onto the first $p$ coordinates:
\begin{equation*}
B^p_{\ell_1} = \left\{\bx \in \R^p ~|~ \exists \bz \in \R^p \, ~\mathrm{s.t.} ~ -\bz_i \leq \bx_i \leq \bz_i \, ~ \forall i, ~ \sum_i \bz_i \leq 1 \right\}.
\end{equation*}
Note that in $\R^{2p}$ with the additional variables $\bz$, we have only $2p+1$ inequalities.
\end{example}

Another example of a polytope that requires many inequalities in a direct description is the \emph{permutahedron} \cite{Zie1995}---the convex hull of all the permutations of the vector $[1, \dots, p]' \in \R^p$.  In fact the permutahedron requires exponentially many linear inequalities in a direct description, while a lifted representation involves $\mathcal{O}(p \log(p))$ additional variables and about $\mathcal{O}(p \log(p))$ inequalities in the higher-dimensional space \cite{Goe2009}.  We refer the reader to the literature on conic representations for other examples (see \cite{GouPT2012} and the references therein), including lifted semidefinite representations.

\subsection*{Hierarchies of Convex Relaxations}

In many cases of interest, convex sets may not have tractable representations.  Lifted representations in such cases have lifting dimensions that are super-polynomially large in the dimension of the original convex set, and thus the associated numerical techniques lead to intractable computational procedures that have super-polynomial runtime with respect to the dimension of the original set.  A prominent example of a convex set that is difficult to compute is the \emph{cut polytope}:
\begin{equation}
\mathrm{CUT}_{m \times m} = \mathrm{conv}\{\bm \bm' ~ |~ \bm \in \{-1,+1\}^m  \}. \label{eq:cut}
\end{equation}
Rank-one signed matrices and their convex combinations are of interest in collaborative filtering and clustering problems (see the section on time-data tradeoffs).  There is no known tractable representation of the cut polytope---lifted linear or semidefinite representations have lifting dimensions that are super-polynomial in size.  Such computational issues have led to a large literature on approximating intractable convex sets by tractable ones.  For the purposes of this paper, and following the dominant trend in the literature, we focus on outer approximations.  For example, the elliptope \eqref{eq:elliptope} is an outer relaxation of the cut polytope, and it has been employed in approximation algorithms for intractable combinatorial optimization problems such as finding the maximum-weight cut in a graph \cite{GoeW1995}.  More generally, one can imagine a hierarchy of increasingly tighter approximations $\{\bc_i\}$ of a convex set $\bc$ as follows:
\begin{equation*}
\bc \subseteq \cdots \subseteq \bc_3 \subseteq \bc_2 \subseteq \bc_1.
\end{equation*}
There exist several mechanisms for deriving such hierarchies, and we describe three frameworks here.

In the first framework, which was developed by Sherali and Adams \cite{SheA1990}, the set $\bc$ is assumed to be polyhedral and each element of the family $\{\bc_i\}$ is also polyhedral.  Specifically, each $\bc_i$ is expressed via a lifted LP representation.  Tighter approximations are obtained by resorting to larger-sized lifts so that the lifting dimension increases with the level $i$ in the hierarchy.  The second framework is similar in spirit to the first one, but now the set $\bc$ is a convex basic, closed semialgebraic set\footnote{A basic, closed semialgebraic set is the collection of solutions of a system of polynomial equations and polynomial inequalities \cite{BocCR1998}.} and the approximations $\{\bc_i\}$ are given by lifted SDP representations.  Again the lifting dimension increases with the level $i$ in the hierarchy.  This method was initially pioneered by Parrilo \cite{Par2000,Par2003} and by Lasserre \cite{Las2001}, and it was studied in greater detail subsequently by Gouveia et al. \cite{GouPT2010}.  Both these first and second frameworks are similar in spirit in that tighter approximations are obtained by via lifted representations with successively larger lifting dimensions.  The third framework we mention here is qualitatively different from the first two.  Suppose $\bc$ is a convex set that has a $\bk$-representation---by successively weakening the cone $\bk$ itself one obtains increasingly weaker approximations to $\bc$.  Specifically, we consider the setting in which the cone $\bk$ is a hyperbolicity cone \cite{Ren2006}.  Such cones have rich geometric and algebraic structure, and their boundary is given in terms of the vanishing of hyperbolic polynomials.  They include the orthant and the semidefinite cone as special cases.  We do not go into further technical details and formal definitions of these cones here, and instead refer the interested reader to \cite{Ren2006}.  The main idea is that one can obtain a family of relaxations $\{\bk_i\}$ to a hyperbolicity cone $\bk \subseteq \R^p$ where each $\bk_i$ is a convex cone (in fact, hyperbolic) and is a subset of $\R^p$:
\begin{equation}
\bk \subseteq \cdots \subseteq \bk_3 \subseteq \bk_2 \subseteq \bk_1. \label{eq:conerelax}
\end{equation}
These outer conic approximations are obtained by taking certain derivatives of the hyperbolic polynomial used to define the original cone $\bk$---see \cite{Ren2006} for more details.  One then constructs a hierarchy of approximations $\{\bc_i\}$ to $\bc$ by replacing the cone $\bk$ in the representation of $\bc$ by the family of conic approximations $\{\bk_i\}$.  From \eqref{eq:cone-rep} and \eqref{eq:conerelax} it is clear that the $\{\bc_i\}$ so defined satisfy $\bc \subseteq \cdots \subseteq \bc_3 \subseteq \bc_2 \subseteq \bc_1$.

The important point in these three frameworks is that the family of approximations $\{\bc_i\}$ obtained in each case is ordered both by approximation quality as well as by computational complexity; that is, the weaker approximations in the hierarchy are also the ones that are more tractable to compute.  This observation leads to an algorithm weakening mechanism that is useful for processing larger datasets more coarsely.  As demonstrated concretely in the next section, the estimator \eqref{eq:shrink} based on a weaker approximation to $\bc$ can provide the same statistical performance as one based on a stronger approximation to $\bc$ provided the former estimator is evaluated with more data.  The upshot is that the first estimator is more tractable to compute than the second.  Thus, we obtain a technique for reducing the runtime required to process a larger dataset.

\section*{Estimation via Convex Optimization}

In this section we investigate the statistical properties of the estimator \eqref{eq:shrink} for the denoising problem \eqref{eq:denoise}.  The signal set $\bs$ in \eqref{eq:denoise} differs based on the application of interest.  For example $\bs$ may be the set of sparse vectors in a fixed basis, which could correspond to the problem of denoising sparse vectors in wavelet bases \cite{Don1995}.  The signal set $\bs$ may be the set of low-rank matrices, which leads to problems of collaborative filtering \cite{SreS2005}.  Finally, $\bs$ may be a set of permutation matrices corresponding to rankings over a collection of items.  Our analysis in this section is general, and is applicable to these and other settings (see the section on time-data tradeoffs for concrete examples).  In some denoising problems, one is interested in noise models other than Gaussian.  We comment on the performance of the estimator \eqref{eq:shrink} in settings with non-Gaussian noise, although we primarily focus on the Gaussian case for simplicity.

\subsection*{Convex Programming Estimators}

In order to analyze the performance of the estimator \eqref{eq:shrink}, we introduce a few concepts from convex analysis \cite{Roc1996}.  Given a closed convex set $\bc \in \R^p$ and a point $\ba \in \bc$ we define the \emph{tangent cone} at $\ba$ with respect to $\bc$ as
\begin{equation}
T_{\bc}(\ba) = \mathrm{cone}\{\bb-\ba~|~ \bb \in \bc\}. \label{eq:tcone}
\end{equation}
Here $\mathrm{cone}(\cdot)$ refers to the conic hull of a set obtained by taking nonnegative linear combinations of elements of the set.  The cone $T_{\bc}(\ba)$ is the set of directions to points in $\bc$ from the point $\ba$.  The \emph{polar} $\bk^\ast \subseteq \R^p$ of a cone $\bk \subseteq \R^p$ is the cone
\begin{equation*}
\bk^\ast = \{\bh \in \R^p ~|~ \langle \bh,\bd \rangle \leq 0 \; \forall \bd \in \bk\}.
\end{equation*}
The \emph{normal cone} $N_{\bc}(\ba)$ at $\ba$ with respect to the convex set $\bc$ is the polar cone of the tangent cone $T_{\bc}(\ba)$:
\begin{equation}
N_{\bc}(\ba) = T_{\bc}(\ba)^\ast. \label{eq:ncone}
\end{equation}
Thus, the normal cone consists of vectors that form an obtuse angle with every vector in the tangent cone $T_{\mathcal{C}}(\bx)$.  Both the tangent and normal cones are convex cones.

A key quantity that will appear in our error bounds is the following notion of the ``complexity'' or ``size'' of a tangent cone:
\begin{definition}
The \emph{Gaussian squared-complexity} of a set $\mathcal{D} \in \R^p$ is defined as:
\begin{equation*}
\G(\mathcal{D}) = \mathbb{E}\left[\sup_{\ba \in \mathcal{D}} ~ \langle \ba, \bg \rangle ^2 \right],
\end{equation*}
where the expectation is with respect to $\bg \sim \mathcal{N}(0,I_{p \times p})$.
\end{definition}
This quantity is closely related to the Gaussian complexity of a set \cite{Dud1967,BarM2002} which consists of no squaring of the term inside the expectation.  The Gaussian squared-complexity shares many properties in common with the Gaussian complexity, and we describe those that are relevant to this paper in the next subsection.  Specifically, we discuss methods to estimate this quantity for sets $\mathcal{D}$ that have some structure.

With these definitions and letting $B_{\ell_2}^p$ denote the $\ell_2$ ball in $\R^p$, we have the following result on the error between $\hat{\bx}_n(C)$ and $\bxa$:
\begin{proposition}\label{prop:denoise}
For $\bxa \in \bs \subset \R^p$ and with $\bc \subseteq \R^p$ convex such that $\bs \subseteq \bc$, we have the error bound
\begin{equation*}
\mathbb{E}\left[\|\bxa - \hat{\bx}_n(\bc)\|_{\ell_2}^2 \right] \leq \tfrac{\sigma^2}{n} \G(T_\bc(\bxa) \cap B_{\ell_2}^p).
\end{equation*}
\end{proposition}

\noindent \textbf{Proof}:  We have that $\bar{\by} = \bxa + \tfrac{\sigma}{\sqrt{n}} \bz$.  We begin by establishing a bound that is derived by conditioning on $\bz = \tz$.  Subsequently, taking expectations concludes the proof.  We have from the optimality conditions \cite{Roc1996} of the convex program \eqref{eq:shrink} that
\begin{equation*}
\bxa + \tfrac{\sigma}{n} \tz - \hat{\bx}_n(\bc) |_{\bz = \tz} \in N_\bc(\hat{\bx}_n(\bc) |_{\bz = \tz}).
\end{equation*}
Here $\hat{\bx}_n(\bc) |_{\bz = \tz}$ represents the optimal value of \eqref{eq:shrink} conditioned on $\bz = \tz$.  As $\bxa \in \bs \subseteq \bc$, we have that $\bxa - \hat{\bx}_n(\bc) |_{ \bz = \tz} \in T_\bc(\hat{\bx}_n(\bc) |_{ \bz = \tz})$.  Since the normal and tangent cones are polar to each other, we have that
\begin{equation*}
\langle \bxa + \tfrac{\sigma}{\sqrt{n}} \tz - \hat{\bx}_n(\bc) |_{\bz = \tz}, \bxa - \hat{\bx}_n(\bc) |_{ \bz = \tz} \rangle \leq 0.
\end{equation*}
It then follows that
\begin{eqnarray*}
\|\bxa &-& \hat{\bx}_n(\bc) |_{\bz = \tz}\|_{\ell_2}^2 \\ &\leq& \tfrac{\sigma}{\sqrt{n}} \left\langle \hat{\bx}_n(\bc) |_{ \bz = \tz} - \bxa, \tz \right\rangle \\ &=& \tfrac{\sigma}{\sqrt{n}} \|\hat{\bx}_n(\bc) |_{ \bz = \tz} - \bxa\|_{\ell_2} \left\langle \frac{\hat{\bx}_n(\bc) |_{ \bz = \tz} - \bxa}{\|\hat{\bx}_n(\bc) |_{ \bz = \tz} - \bxa\|_{\ell_2}}, \tz \right\rangle \\ &\leq& \tfrac{\sigma}{\sqrt{n}} \|\hat{\bx}_n(\bc) |_{ \bz = \tz} - \bxa\|_{\ell_2} ~ \left[\sup_{\bd \in T_\bc(\bxa), \|\bd\|_{\ell_2} \leq 1} ~ \langle \bd, \tz \rangle \right].
\end{eqnarray*}
Dividing both sides by $\|\hat{\bx}_n(\bc) |_{ \bz = \tz} - \bxa\|_{\ell_2}$, then squaring both sides, and finally taking expectations completes the proof. $\square$

Note that the basic structure of the error bound provided by the estimator \eqref{eq:shrink} in fact holds for an \emph{arbitrary} distribution on the noise $\bz$ with the Gaussian squared-complexity suitably modified.  However, we focus for the rest of this paper on the Gaussian case, $\bz \sim \mathcal{N}(0, I_{p \times p})$.

To summarize in words, the mean squared error is bounded by the noise variance times the Gaussian squared-complexity of the normalized tangent cone with respect to $\bc$ at the true parameter $\bxa$.  Essentially it measures the amount of noise restricted to the tangent cone, which is intuitively reasonable as only the noise that moves one away from $\bxa$ in a feasible direction in $\bc$ must contribute towards the error.  Therefore, if the convex constraint set $\bc$ is ``sharp'' at $\bxa$ so that the cone $T_\bc(\bxa)$ is ``narrow,'' then the error is small.  At the other extreme, if the constraint set $\bc = \R^p$ then the error is $\tfrac{\sigma^2}{n} p$ as one would expect.

While Proposition~\ref{prop:denoise} is useful and indeed it will suffice for the purposes of demonstrating time-data tradeoffs in the next section, there are a couple of shortcomings in the result as stated.  First, suppose the signal set $\bs$ is contained in a ball around the origin with the radius of the ball being small relative to the noise variance $\tfrac{\sigma^2 }{n}$.  In such a setting, the estimator $\hat{\bx} = 0$ leads to a smaller mean squared error than one would obtain from Proposition~\ref{prop:denoise}.  Second, and somewhat more subtly, suppose that one doesn't have a perfect bound on the size of the signal set.  For concreteness, consider a setting in which $\bs$ is a set of sparse vectors with bounded $\ell_1$ norm, in which case a good choice for the constraint set $\bc$ in the estimator \eqref{eq:denoise} is an appropriately scaled $\ell_1$ ball.  However, if we do not know the $\ell_1$ norm of $\bxa$ \emph{a priori}, then we may end up employing a constraint set $\bc$ such that $\bxa$ does not belong to $\bc$ (so $\bxa$ is an infeasible solution) or that $\bxa$ lies strictly in the interior of $\bc$ (hence, $T_\bc(\bxa)$ is all of $\R^p$).  Both these situations are undesirable as they limit the applicability of Proposition~\ref{prop:denoise} and provide very loose error bounds.  The following result addresses these shortcomings by weakening the assumptions of Proposition~\ref{prop:denoise}:
\begin{proposition} \label{prop:denoise2}
Let $\bxa \in \bs \subset \R^p$ and let $\bc \subseteq \R^p$.  Suppose there exists a point $\tilde{\bx} \in \bc$ such that
$\bc - \tilde{\bx} = Q_1 \oplus Q_2$, with $Q_1, Q_2$ lying in orthogonal subspaces of $\R^p$ and $Q_2 \subseteq \alpha B^p_{\ell_2}$ for $\alpha \geq 0$.  Then we have that
\begin{eqnarray*}
\mathbb{E}\left[\|\bxa - \hat{\bx}_n(\bc)\|_{\ell_2}^2 \right] \leq 6 \Big[\tfrac{\sigma^2}{n} \G(\mathrm{cone}(Q_1) \cap B^p_{\ell_2}) + \\ \|\bxa - \tilde{\bx}\|_{\ell_2}^2 + \alpha^2 \Big].
\end{eqnarray*}
Here $\mathrm{cone}(Q_1)$ is the conic hull of $Q_1$.
\end{proposition}
The proof of this result is presented in the Supplementary Information.  A number of remarks are in order here.  With respect to the first shortcoming in Proposition~\ref{prop:denoise} stated above, if $\bc$ is chosen such that $\bs \subset \bc$ one can set $\tilde{\bx}=\bxa, Q_1 = 0$, and $Q_2 = \bc - \tilde{\bx}$ in Proposition~\ref{prop:denoise2}, and readily obtain a bound that scales only with the diameter of the convex constraint set $\bc$.  With regard to the second shortcoming in Proposition~\ref{prop:denoise} described above, if a point $\tilde{\bx} \in \bc$ near $\bxa$ has a narrow tangent cone $T_\bc(\tilde{\bx})$, then one can provide an error bound with respect to $\G(T_\bc(\tilde{\bx}))$ with an extra additive term that depends on $\|\bxa - \tilde{\bx}\|_{\ell_2}^2$---this is done by setting $Q_1 = \bc-\tilde{\bx}$ and $Q_2 = 0$ (thus, $\alpha = 0$) in Proposition~\ref{prop:denoise2}.  More generally, Proposition~\ref{prop:denoise2} incorporates both these improvements in a single error bound with respect to an arbitrary point $\tilde{\bx} \in \bc$; thus, one can further optimize the error bound over the choice of $\tilde{\bx} \in \bc$ (as well as the choice of the decomposition $Q_1$ and $Q_2$).

\subsection*{Properties and Computation of Gaussian Squared-Complexity}

We record some properties of the Gaussian squared-complexity that are subsequently useful when we demonstrate concrete time-data tradeoffs.  It is clear that $\G(\cdot)$ is monotonic with respect to set nesting, i.e., $\G(\mathcal{D}_1) \leq \G(\mathcal{D}_2)$ for sets $\mathcal{D}_1 \subseteq \mathcal{D}_2$.  If $\mathcal{D}$ is a subspace then one can check that $\G(\mathcal{D}) = \mathrm{dim}(\mathcal{D})$.  In order to estimate squared-complexities of families of cones, one can imagine appealing to techniques similar to those used for estimating Gaussian complexities of sets \cite{Dud1967,BarM2002}.  Most prominent among these are arguments based on covering number and metric entropy bounds.  However, these arguments are frequently not sharp and introduce extraneous log-factors in the resulting error bounds.

In a recent paper \cite{ChaRPW2012}, sharp upper bounds on the Gaussian complexities of
normalized cones have been established for families of cones of interest in a class of
linear inverse problems.  The (square of the) Gaussian complexity can be upper bounded
by the Gaussian squared-complexity $\G(\mathcal{D})$ via Jensen's inequality:
\begin{equation*}
\mathbb{E} \left[\sup_{\bd \in \mathcal{D}} ~ \langle \bd, \bg \rangle\right]^2 \leq \G(\mathcal{D}),
\end{equation*}
where $\bg$ is a standard normal vector.  In fact most of these bounds in \cite{ChaRPW2012}
were obtained by bounding $\G(\mathcal{D})$ and thus they are directly relevant to our
setting.  In the rest of this section we present the bounds on $\G(\mathcal{D})$
from \cite{ChaRPW2012} that will be used in this paper, deferring to that paper
for proofs in most cases.  In some cases the proofs do require modifications with
respect to their counterparts in \cite{ChaRPW2012}, and for these cases we give full
proofs in the Supplementary Information.

The first result is a direct consequence of convex duality and provides a fruitful general technique to compute sharp estimates of Gaussian squared-complexities.  Let $\mathrm{dist}(\ba,\mathcal{D})$ denote the $\ell_2$ distance from a point $\ba$ to the set $\mathcal{D}$.
\begin{lemma}\cite{ChaRPW2012} \label{lemm:dual}
Let $\bk \subseteq \R^p$ be a convex cone and let $\bk^\ast \subseteq \R^p$ be its polar.  Then we have for any $\ba \in \R^p$ that
\begin{equation*}
\sup_{\bd \in \bk \cap B^p_{\ell_2}} ~ \langle \bd, \ba \rangle = \mathrm{dist}(\ba,\bk^\ast).
\end{equation*}
\end{lemma}
Therefore we have the following result as a simple corollary.
\begin{corollary} \label{cor:gsc2dist}
Let $\bk \subseteq \R^p$ be a convex cone and let $\bk^\ast \subseteq \R^p$ be its polar.  For $\bg \sim \mathcal{N}(0,I_{p\times p})$ we have that
\begin{equation*}
\G(\bk \cap B^p_{\ell_2}) = \mathbb{E}\left[\mathrm{dist}(\bg,\bk^\ast)^2 \right].
\end{equation*}
\end{corollary}

Based on the duality result of Lemma~\ref{lemm:dual} and Corollary~\ref{cor:gsc2dist}, one can compute the following sharp bounds on the Gaussian squared-complexities of tangent cones with respect to the $\ell_1$ norm and nuclear norm balls.  These are especially relevant when one wishes to estimate sparse signals or low-rank matrices---in these settings the $\ell_1$ and nuclear norm balls serve as useful constraint sets for denoising as the tangent cones with respect to these sets at sparse vectors and at low-rank matrices are particularly narrow.  Both these results are used when we describe time-data tradeoffs.
\begin{proposition}\cite{ChaRPW2012} \label{prop:l1}
Let $\bx \in \R^p$ be a vector containing $s$ nonzero entries.  Let $T$ be the tangent cone at $\bx$ with respect to an $\ell_1$ norm ball scaled so that $\bx$ lies on the boundary of the ball, i.e., a scaling of the unit $\ell_1$ norm ball by a factor $\|\bx\|_{\ell_1}$.  Then
\begin{equation*}
\G(T\cap B^p_{\ell_2}) \leq 2s\log(\tfrac{p}{s}) + \tfrac{5}{4} s.
\end{equation*}
\end{proposition}
Next we state a result for low-rank matrices and the nuclear norm ball.
\begin{proposition}\cite{ChaRPW2012} \label{prop:nuclear}
Let $\mathbf{X} \in \R^{m_1 \times m_2}$ be matrix of rank $r$.  Let $T$ be the tangent cone at $\mathbf{X}$ with respect to a nuclear norm ball scaled so that $\mathbf{X}$ lies on the boundary of the ball, i.e., a scaling of the unit nuclear norm ball by a factor equal to the nuclear norm of $\mathbf{X}$.  Then
\begin{equation*}
\G(T\cap B^{m_1 m_2}_{\ell_2}) \leq 3r(m_1+m_2-r).
\end{equation*}
\end{proposition}

Next we state and prove a result that allows us to estimate Gaussian squared-complexities of general cones.  The bound is based on the volume of the dual of the cone of interest, and the proof involves an appeal to Gaussian isoperimetry \cite{Led2000}.  A similar result on Gaussian complexities of cones (without the square) was proved in \cite{ChaRPW2012}, but that result does not directly imply our statement and we therefore give a complete self-contained proof in the Supplementary Information.  The volume of a cone is assumed to be normalized (between zero and one) so we consider the relative fraction of a unit Euclidean sphere that is covered by a cone:
\begin{proposition}\label{prop:isop}
Let $\bk \subset \R^p$ be a cone such that its polar $\bk^\ast \subset \R^p$ has a normalized volume of $\mu \in (\tfrac{1}{4} \exp\{-p/20\}, \tfrac{1}{4e^2})$.  For $p \geq 12$, we have that
\begin{equation*}
\G(\bk \cap B^p_{\ell_2}) \leq 20 \log\left(\tfrac{1}{4\mu}\right).
\end{equation*}
\end{proposition}
If a cone is narrow then its polar will be wide leading to a large value of $\mu$ and hence a small quantity of the right-hand-side of the bound.  This result leads to bounds on Gaussian squared-complexity in settings in which one can easily obtain estimates of volumes.  One setting in which such estimates are easily obtained is the case of tangent cones with respect to vertex transitive polytopes.  We recall that a vertex transitive polytope \cite{Zie1995} is one in which there exists a symmetry of the polytope for each pair of vertices mapping the two vertices isomorphically to each other. Roughly speaking, all the vertices in such polytopes are the same.  Some examples include the cross-polytope (the $\ell_1$ norm ball), the simplex \eqref{eq:simplex}, the hypercube (the $\ell_\infty$ norm ball), and many polytopes generated by the action of groups \cite{SanSS2009}.  We will see many examples of such polytopes in our examples on time-data tradeoffs, and thus we will appeal to the following corollary repeatedly:
\begin{corollary}\label{cor:sym}
Suppose that $\mathcal{P} \in \R^p$ is a vertex transitive polytope with $v$ vertices and let $\bx$ be a vertex of this polytope.  If $4e^2 \leq v \leq 4\exp\{p/20\}$ then
\begin{equation*}
\G(T_{\mathcal{P}}(\bx)) \leq 20 \log(v/4).
\end{equation*}
\end{corollary}
\noindent \textbf{Proof}: The normal cones at the vertices of $\mathcal{P}$ partition $\R^p$.  If the polytope is vertex-transitive, then the normal cones are all equivalent to each other (up to orthogonal transformations).  Consequently, the (normalized) volume of the normal cone at any vertex is $1 / v$.  Since the normal cone at a vertex is polar to the tangent cone, we have the desired result from Proposition~\ref{prop:isop}. $\square$

\begin{figure}
\begin{center}
\epsfig{file=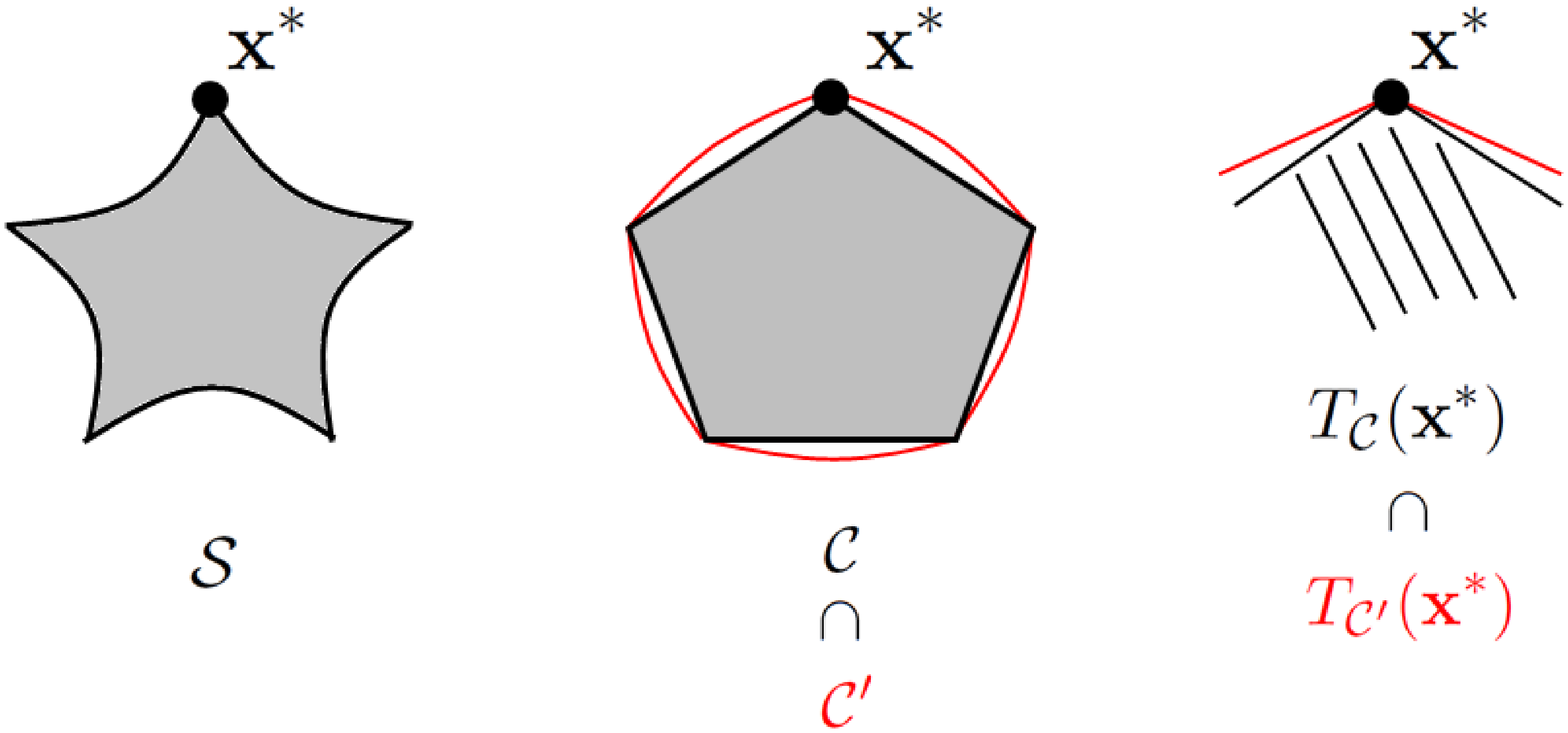,width=8cm,height=4cm} \caption{(left) A signal set $\bs$ consisting of $\bxa$; (middle) Two convex constraint sets $\bc$ and $\bc'$, where $\bc$ is the convex hull of $\bs$ and $\bc'$ is a relaxation that is more efficiently computable than $\bc$; (right) The tangent cone $T_\bc(\bxa)$ is contained inside the tangent cone $T_{\bc'}(\bxa)$.  Consequently, the Gaussian squared-complexity $\G(T_\bc(\bxa) \cap B_{\ell_2}^p)$ is smaller than the complexity $\G(T_{\bc'}(\bxa) \cap B_{\ell_2}^p)$, so that the estimator $\hat{\bx}_n(\bc)$ requires fewer samples than the estimator $\hat{\bx}_n(\bc')$ for a risk of at most $1$.} \label{fig:relax}
\end{center}
\end{figure}

\section*{Time-Data Tradeoffs}

\subsection*{Preliminaries}

We now turn our attention to giving examples of time-data tradeoffs in denoising problems via convex relaxation.  As described previously, we must set a desired risk in order to realize a time-data tradeoff---in the examples in the rest of this section, we will fix the desired risk to be equal to $1$ independent of the problem dimension $p$.  Thus, these denoising problems lie in the time-data complexity classes $\td(t(p),n(p),1)$ for different runtime constraints $t(p)$ and sample budgets $n(p)$.  The following corollary gives the number of samples required to obtain a mean squared error of $1$ via convex optimization in our denoising setup:

\begin{corollary}\label{cor:risk}
For $\bxa \in \bs$ and with $\bs \subseteq \bc$, if
\begin{equation*}
n \geq \sigma^2 \G(T_\bc(\bxa) \cap B^p_{\ell_2}),
\end{equation*}
then $\mathbb{E}\left[\|\bxa - \hat{\bx}_n(\bc)\|_{\ell_2}^2 \right] \leq 1$.
\end{corollary}

\noindent \textbf{Proof}: The result follows by a rearrangement of the terms in the bound in Proposition~\ref{prop:denoise}. $\square$

This corollary states that if we have access to a dataset with $n$ samples, then we can use any convex constraint set $\bc$ such that the term on the right-hand-side in the corollary is smaller than $n$.  Recalling that larger constraint sets $\bc$ lead to larger tangent cones $T_\bc$, we observe that if $n$ is large one can potentially use very weak (and computationally inexpensive) relaxations and still obtain a risk of $1$.  This observation, combined with the important point that the hierarchies of convex relaxations described previously are simultaneously ordered both by approximation quality and by computational tractability, allows us to realize a time-data tradeoff by using convex relaxation as an algorithm weakening mechanism.  See Figure~2 for a simple demonstration.

We further consider settings with $\sigma^2 = 1$ and in which our signal sets $\bs \subseteq \R^p$ consist of elements that have Euclidean norm on the order of $\sqrt{p}$ (measured from the centroid of $\bs$).  In such regimes the James-Stein shrinkage estimator \cite{JamS1961} offers about the same level of performance as the maximum-likelihood estimator, and both these are outperformed in statistical risk by nonlinear estimators of the form \eqref{eq:shrink} based on convex optimization.

Finally, we briefly remark on the runtimes of our estimators.  The runtime for each of the procedures below is calculated by adding the number of operations required to compute the sample mean $\bar{\by}$ and the number of operations to solve \eqref{eq:shrink} to some accuracy.  Hence if the number of samples used is $n$ and if $f_\bc(p)$ denotes the number of operations required to project $\bar{\by}$ onto $\bc$, then the total runtime is $np + f_\bc(p)$.  Thus the number of samples enters the runtime calculations as just an additive term. As we process larger datasets, the first term in this calculation becomes larger but this increase is offset by a more substantial decrease in the second term due to the use of a computationally tractable convex relaxation.  We note that such a runtime calculation extends to more general inference problems in which one employs estimators of the form \eqref{eq:shrink} but with different loss functions in the objective---specifically, the runtime is calculated as above so long as the loss function depends only on some sufficient statistic computed from the data.  If the loss function is instead of the form $\sum_{i=1}^n \ell(\bx; \by_i)$ and it cannot be summarized via a sufficient statistic of the data $\{\by_i\}_{i=1}^n$, then the number of samples enters the runtime computation in a multiplicative manner in the number of operations $f_\bc(p)$ required to compute the convex programming estimator.

\subsection*{Example 1: Denoising Signed Matrices}

We consider the problem of recovering signed matrices corrupted by noise:
\begin{equation*}
\bs = \{\ba \ba' ~|~ \ba \in \{-1,+1\}^{\sqrt{p}} \}.
\end{equation*}
We have $\ba \in \R^{\sqrt{p}}$ so that $\bs \subseteq \R^p$.  Inferring such signals is of interest in collaborative filtering where one wishes to approximate matrices as the sum of a small number of rank-one signed matrices \cite{SreS2005}.  Such matrices may represent, for example, the movie preferences of users as in the Netflix problem.

The tightest convex constraint that one could employ in this case is $\bc = \mathrm{conv}(\bs)$, which is the cut polytope \eqref{eq:cut}.  In order to obtain a risk of $1$ with this constraint, one requires $n = c_1 \sqrt{p}$ by applying Corollary~\ref{cor:sym} and Corollary~\ref{cor:risk} based on the symmetry of the cut polytope.  The cut polytope is in general intractable to compute.  Hence the best known algorithms to project onto $\bc$ would require runtime that is super-polynomial in $p$.  Consequently, the total runtime of this algorithm is $c_1 p^{1.5} + \mathrm{super}$-$\mathrm{poly}(p)$.

A commonly used tractable relaxation of the cut polytope is the elliptope \eqref{eq:elliptope}.  By computing the Gaussian squared-complexity of the tangent cones at rank-one signed matrices with respect to this set, it is possible to show that $n = c_2 \sqrt{p}$ leads to a risk of $1$ (with $c_2 > c_1$).  Further, interior-point based convex optimization algorithms for solving \eqref{eq:shrink} that exploit the special structure of the elliptope require $\mathcal{O}(p^{2.25})$ operations\footnote{The exponent is a result of the manner in which we define our signal set so that a rank-one signed matrix lives in $\R^p$.} \cite{BoyV2004,BenN2001,Hig2002}.  Hence the total runtime of this procedure is $c_2 p^{1.5} + \mathcal{O}(p^{2.25})$.

Finally, an even weaker relaxation of the cut polytope than the elliptope is the unit ball of the nuclear norm scaled by a factor of $\sqrt{p}$---one can verify that the elements of $\bs$ lie on the boundary of this set, and are in fact extreme points.  Appealing to Proposition~\ref{prop:nuclear} (using the fact that the elements of $\bs$ are rank-one matrices) and Corollary~\ref{cor:risk}, we conclude that $n = c_3 \sqrt{p}$ samples provide a mean-squared error of $1$ (with $c_3 > c_2$).  Projecting onto the scaled nuclear norm ball can be done by computing a singular value decomposition (SVD), and then truncating the sequence of singular values in descending order when their cumulative sum exceeds $\sqrt{p}$ (in effect projecting the vector of singular values onto an $\ell_1$ ball of size $\sqrt{p}$).  This operation requires $\mathcal{O}(p^{1.5})$ operations, and thus the total runtime is $c_3 p^{1.5} + \mathcal{O}(p^{1.5})$.

To summarize, the cut-matrix denoising problem lives in the time-data class $\td(\mathrm{super}$-$\mathrm{poly}(p), \allowbreak c_1 \sqrt{p}, 1)$, in $\td(\mathcal{O}(p^{2.25}), \allowbreak c_2 \sqrt{p}, 1)$, and in $\td(\mathcal{O}(p^{1.5}), c_3 \sqrt{p}, 1)$, with constants $c_1 < c_2 < c_3$.

\subsection*{Example 2: Ordering Variables}

In many data analysis tasks, one is given a collection of variables that are suitably ordered so that the population covariance is banded.  Under such a constraint, thresholding the entries of the empirical covariance matrix based on their distance from the diagonal has been shown to be a powerful method for estimation in the high-dimensional setting \cite{BicL2008}.  However, if an ordering of the variables is not known \emph{a priori}, then one must jointly learn an ordering for the variables and estimate their underlying covariance.  As a stylized version of this variable ordering problem, let $M \in \R^{\sqrt{p} \times \sqrt{p}}$ be a known tridiagonal matrix (with Euclidean norm $\mathcal{O}(\sqrt{p})$) and consider the following signal set:
\begin{equation*}
\bs = \{\Pi M \Pi' \;|\; \Pi \; \mathrm{is \; a \;} \sqrt{p} \times \sqrt{p} \; \mathrm{permutation \; matrix}\}.
\end{equation*}
The matrix $M$ here is to be viewed as a covariance matrix.  Thus, the corresponding denoising problem \eqref{eq:denoise} is that we wish to estimate a covariance matrix in the absence of knowledge of the ordering of the underlying variables.  In a real-world scenario one might wish to consider covariance matrices $M$ that belong to some class of banded matrices and then construct $\bs$ as done here, but we stick with the case of a fixed $M$ for simplicity.  Further, the noise in a practical setting is better modeled as coming from a Wishart distribution---again, we focus on the Gaussian case for simplicity.

The tightest convex constraint set that one could employ in this case is the convex hull of $\bs$, which is in general intractable to compute for arbitrary matrices $M$.  For example, if one were able to compute this set in polynomial time for any tridiagonal matrix $M$, one would be able to solve the intractable longest path problem \cite{GarJ1979} (finding the longest path between any two vertices in a graph) in polynomial time.  With this convex constraint set, we find using Corollary~\ref{cor:sym} and Corollary~\ref{cor:risk} that $n = c_1 \sqrt{p} \log(p)$ samples would lead to a risk of $1$.  This follows from the fact that $\mathrm{conv}(\bs)$ is a vertex-transitive polytope with about $(\sqrt{p})!$ vertices.  Thus, the total runtime is $c_1 p^{1.5} \log(p) + \mathrm{super}$-$\mathrm{poly}(p)$.

An efficiently computable relaxation of $\mathrm{conv}(\bs)$ is a scaled $\ell_1$ ball (scaled by the $\ell_1$ norm of $M$).  Appealing to Proposition~\ref{prop:l1} on tangent cones with respect to the $\ell_1$ ball and to Corollary~\ref{cor:risk}, we find that $n = c_2 \sqrt{p} \log(p)$ samples suffice to provide a risk of $1$.  In applying Proposition~\ref{prop:l1}, we note that $M$ is assumed to be tridiagonal and therefore has $\mathcal{O}(\sqrt{p})$ nonzero entries.  The runtime of this procedure is $c_2 p^{1.5} \log(p) + \mathcal{O}(p \log(p))$.

Thus the variable ordering denoising problem belongs to $\td(\mathrm{super}$-$\mathrm{poly}(p), c_1 \sqrt{p} \log(p), 1)$ and to $\td(\mathcal{O}(p^{1.5} \log(p)), \allowbreak c_2 \sqrt{p} \log(p), 1)$, with constants $c_1 < c_2$.

\subsection*{Example 3: Sparse PCA and Network Activity Identification}

As our third example, we consider sparse PCA in which one wishes to learn from samples a sparse eigenvector that contains most of the energy of a covariance matrix.  As a simplified version of this problem, one can imagine a matrix $M \in \R^{\sqrt{p} \times \sqrt{p}}$ with entries equal to $\sqrt{p} / k$ in the top-left $k \times k$ block and zeros elsewhere (so that the Euclidean norm of $M$ is $\sqrt{p}$), and with $\bs$ defined as:
\begin{equation*}
\bs = \{\Pi M \Pi' \;|\; \Pi \;\mathrm{is \; a \;} \sqrt{p} \times \sqrt{p} \; \mathrm{permutation \; matrix}\}.
\end{equation*}
In addition to sparse PCA, such signal sets are also of interest in identifying activity in noisy networks \cite{KolBRS2011}, as well as in related combinatorial optimization problems such as the planted clique problem \cite{GarJ1979}.  In the sparse PCA context, Amini and Wainwright \cite{AmiW2009} study time-data tradeoffs by investigating the sample complexities of two procedures, a simple one based on thresholding and a more sophisticated one based on semidefinite programming.  Kolar et al. \cite{KolBRS2011} investigate the sample complexities of a number of procedures ranging from a combinatorial search method, thresholding, and sparse SVD.  We note that the time-data tradeoffs studied in these two papers \cite{AmiW2009,KolBRS2011} relate to the problem of learning the support of the leading sparse eigenvector; in contrast in our setup the objective is to simply denoise an element of $\bs$.  Further while the Gaussian noise setting is of interest in some of these domains, in a more realistic sparse PCA problem (such as the one considered in \cite{AmiW2009}) the noise is Wishart rather than Gaussian as considered here.  Nevertheless, we stick with our stylized problem setting as it provides some useful insights on time-data tradeoffs.  Finally, the size of the block $k \in \{1,\dots,\sqrt{p}\}$ depends on the application of interest and it is typically far from the extremes $1$ and $\sqrt{p}$---we will consider the case $k \sim p^{1/4}$ for concreteness.\footnote{This setting is an interesting threshold case in the planted clique context \cite{AloKS1998,FeiK2000,AmeV2011} where $k = p^{1/4}$ is the square-root of the number of nodes $\sqrt{p}$ of the graph represented by $M$ (viewed as an adjacency matrix).}

As usual, the tightest convex constraint set one can employ in this setting is the convex hull of $\bs$, which is in general intractable to compute---an efficient characterization of this polytope would lead to an efficient solution of the intractable planted-clique problem (finding a fully connected subgraph inside a larger graph).  Using this convex constraint set gives an estimator that requires about $n = \mathcal{O}(p^{1/4} \log(p))$ samples in order to produce a risk-$1$ estimate.  We obtain this threshold by appealing to Corollary~\ref{cor:sym} and to Corollary~\ref{cor:risk}, and the observation that $\mathrm{conv}(\bs)$ is a vertex-transitive polytope with about ${\sqrt{p} \choose p^{1/4}}$ vertices.  Thus, the overall runtime is $\mathcal{O}( p^{5/4} \log(p)) + \mathrm{super}$-$\mathrm{poly}(p)$.

A convex relaxation of $\mathrm{conv}(\bs)$ is the nuclear norm ball scaled by a factor of $\sqrt{p}$ so that the elements of $\bs$ lie on the boundary.  From Proposition~\ref{prop:nuclear} (observing that the elements of $\bs$ are rank-one matrices) and Corollary~\ref{cor:risk}, we have that $n = c \sqrt{p}$ samples give a risk-$1$ estimate with this procedure.  As computed in the example with cut matrices, the overall runtime of this nuclear norm procedure is $c p^{1.5} + \mathcal{O}(p^{1.5})$.

In conclusion the denoising version of sparse PCA lies in $\td(\mathrm{super}$-$\mathrm{poly}(p), \mathcal{O}(p^{1/4} \log(p)), 1)$ and in $\td(\mathcal{O}(p^{1.5}), \allowbreak \mathcal{O}(\sqrt{p}), \allowbreak 1)$.

\subsection*{Example 4: Estimating Matchings}

As our final example, we consider signals that represent the set of all perfect matchings in the complete graph. A matching is any subset of edges of a graph such that no node of the graph is incident to more than one edge in the subset, and a perfect matching is a subset of edges in which every node is incident to exactly one edge in the subset.  Graph matchings arise in a range of inference problems such as in chemical structure analysis \cite{RouB1979} and in network monitoring \cite{ShoKR1999}.  Letting $M$ be the adjacency matrix of some perfect matching in the complete graph on $\sqrt{p}$ nodes, our signal set in this case is defined as follows:
\begin{equation*}
\bs = p^{1/4} \{\Pi M \Pi' \;|\; \Pi \;\mathrm{is \; a \;} \sqrt{p} \times \sqrt{p} \; \mathrm{permutation \; matrix}\}.
\end{equation*}
The scaling of $p^{1/4}$ ensures that the elements of $\bs$ have Euclidean norm of $\sqrt{p}$.  Note that $\bs \subset \R^p$.  The number of elements in $\bs$ is $\frac{(\sqrt{p})!}{\left(\tfrac{\sqrt{p}}{2}\right)! ~ 2^{\sqrt{p} / 2}}$ when $\sqrt{p}$ is an even number (this number is obtained by computing the product of all the odd integers up to $\sqrt{p}$).

The tightest convex relaxation in this case is the convex hull of $\bs$.  Unlike the previous three cases, projecting onto this convex set is in fact a polynomial-time operation,\footnote{Edmonds' blossom algorithm \cite{Edm1965} for computing maximum-weight matchings in polynomial time leads to a separation oracle for this perfect matching polytope.  Subsequently, Padberg and Rao \cite{PadR1982} developed a faster separation oracle for the perfect matching polytope.  These separation oracles in turn lead to polynomial-time projection algorithms via the ellipsoid method \cite{BenN2001}.} with runtime about $\mathcal{O}(p^5)$.  Appealing to Corollary~\ref{cor:sym}, to Corollary~\ref{cor:risk}, and to the fact that $\mathrm{conv}(\bs)$ is a vertex-transitive polytope, we have that $n=c_1 \sqrt{p} \log(p)$ samples provides a risk-$1$ estimate.  Hence the overall runtime is $c_1 p^{1.5} \log(p) + \mathcal{O}(p^5)$.

A tractable relaxation of the perfect matching polytope is a hypersimplex \cite{Zie1995}, obtained by taking the convex hull of all $\sqrt{p} \times \sqrt{p}$ matrices consisting of $\sqrt{p}$ ones and the other entries being equal to zero.  We scale this hypersimplex by a factor of $p^{1/4}$ so that the elements of $\bs$ are on the boundary.  The hypersimplex is also a vertex-transitive polytope like the perfect matching polytope, but with about ${p \choose \sqrt{p}}$ entries.  Hence from Corollary~\ref{cor:sym} and from Corollary~\ref{cor:risk}, we have that $n=c_2 \sqrt{p} \log(p)$ samples will provide a risk-$1$ estimate.  Further, projecting onto the hypersimplex is a very efficient operation based on sorting and has a runtime of $\mathcal{O}(p \log(p))$.  Consequently the total runtime of this procedure is $c_2 p^{1.5} \log(p) + \mathcal{O}(p \log(p))$.

In summary, the matching estimation problem is a member of $\td(\mathcal{O}(p^5) , c_1 \sqrt{p} \log(p), 1)$ and of $\td(\mathcal{O}(p^{1.5} \log(p)) , \allowbreak c_2 \sqrt{p} \log(p), \allowbreak 1)$ with constants $c_1 < c_2$.

\subsection*{Some Observations}

A curious observation that we may take away from these examples is that it is possible to obtain substantial speedups computationally with just a constant factor increase in the size of the dataset.  This suggests that in settings in which obtaining additional data is inexpensive, it may be more economical to procure more data and employ a more basic computational infrastructure rather than to process limited data using powerful and expensive computers.

Our second observation is relevant to all the examples above but we highlight it in the context of denoising cut matrices.  In that setting one can use an even weaker relaxation of the cut polytope than the nuclear norm ball, such as the Euclidean ball (suitably scaled).  While projection onto this set is extremely efficient (requiring $\mathcal{O}(p)$ operations as opposed to $\mathcal{O}(p^{1.5})$ operations for projecting onto the nuclear norm ball), the number of samples required to achieve a risk of $1$ with this approach is $\mathcal{O}(p)$---computing the sample mean with so many samples requires $\mathcal{O}(p^2)$ operations, which leads to an overall runtime that is greater than the runtime $\mathcal{O}(p^{1.5})$ for the nuclear norm approach. This point highlights an important tradeoff---if our choice of algorithms is between nuclear norm projection and Euclidean projection, and if we are in fact given access to $\mathcal{O}(p)$ data samples, it makes sense computationally to retain only $\mathcal{O}(\sqrt{p})$ samples for the nuclear norm procedure and throw away the remaining data.  This provides a concrete illustration of several key issues.  Aggregating massive datasets can frequently be very expensive computationally (relative to the other subsequent processing), and the number of operations required for this step must be taken into account.\footnote{Note that the aggregation step is more time-consuming than the subsequent projection step in the $\ell_1$-ball projection procedure for ordering variables and in the hypersimplex projection method for denoising matchings.}  Consequently, in some cases it may make sense to throw away some data if preprocessing the full massive dataset is time-consuming.  Hence, one may not be able to avail oneself of weaker post-aggregation algorithms if these methods require such a large amount of data to achieve a desired risk that the aggregation step is expensive.  This point goes back to the ``floor'' in Figure~1 in which one imagines a cutoff in the number of samples beyond which more data are not helpful in reducing computational runtime.  Such a threshold, of course, depends on the space of algorithms one employs, and in the cut polytope context with the particular algorithms considered here the threshold occurs at $\mathcal{O}(\sqrt{p})$ samples.

\section*{Conclusions}

In this paper we considered the problem of reducing the computational complexity of an inference task as one has access to larger datasets.  The traditional goal in the theory of statistical inference is to understand the tradeoff in an estimation problem between the amount of data available and the risk attainable via some class of procedures.  In an age of plentiful data in many settings and computational resources being the principal bottleneck, we believe that an increasingly important objective is to investigate the tradeoffs between computational and sample complexities.  As one pursues this line of thinking, it becomes clear that a central theme must be the ability to weaken an inference procedure as one has access to larger datasets.  Accordingly, we proposed convex relaxation as an algorithm weakening mechanism, and we investigated its efficacy in a class of denoising tasks.  Our results suggest that such methods are especially effective in achieving time-data tradeoffs in high-dimensional parameter estimation.

We close our discussion by outlining some exciting future research directions.  As algorithm weakening is central to the viewpoint described in this paper, it should come as no surprise that several of the directions listed below involve interaction with important themes in computer science.

\noindent \paragraph{Computation with streaming data}  In many massive data problems, one is presented with a stream of input data rather than a large fixed dataset, and an estimate may be desired after a fixed amount of time independent of the rate of the input stream.  In such a setting an alternative viewpoint to the one presented in this paper might be more appropriate.  Specifically, rather than keeping the risk fixed, one would keep the runtime fixed and trade off the risk with the rate of the input stream.  One can imagine algorithm weakening mechanisms, dependent on the rate of the data stream, in which the initial data points are processed using sophisticated algorithms and subsequent samples are processed more coarsely.  Understanding the tradeoffs in such a setting is of interest in a range of applications.

\noindent \paragraph{Alternative algorithm weakening mechanisms} The notion of weakening an inference algorithm is key to realizing a time-data tradeoff.  While convex relaxation methods provide a powerful and general approach, a number of other weakening mechanisms are potentially relevant.  For example, processing data more coarsely by quantization, dimension reduction, and clustering may be natural in some contexts.  Coresets, which originated in the computational geometry community, summarize a large set of points via a small collection (see, for example, \cite{FelL2011} and the references therein), and they could also provide a powerful algorithm weakening mechanism.  Finally, we would like to mention a computer hardware concept that has implications for massive data analysis.  A recent approach to designing computer chips is premised on the idea that many tasks do not require extremely accurate computation \cite{Bat}.  If one is willing to tolerate small, random errors in arithmetic computations (e.g., addition, multiplication), it may be possible to design chips that consume less power and are faster than traditional, more accurate chips.  Translated to a data analysis context, such design principles may provide a hardware-based algorithm weakening mechanism.

\noindent \paragraph{Measuring quality of approximation of convex sets} In the mathematical optimization and theoretical computer science communities, relaxations of convex sets have provided a powerful toolbox for designing approximation algorithms for intractable problems, most notably those arising in combinatorial optimization.  The manner in which the quality of a relaxation translates to the quality of an approximation algorithm is usually quantified based on the \emph{integrality gap} between the original convex set and its approximation \cite{Vaz2004}.  However, the quantity of interest in a statistical inference context in characterizing the quality of approximations is based on ratios of Gaussian squared-complexities of tangent cones.  These two quantifications can be radically different---indeed, several of the relaxations presented in our time-data tradeoff examples that are useful in an inferential setting would provide poor performance in a combinatorial optimization context.  More broadly, those examples demonstrate that weak relaxations frequently provide as good estimation performance as tighter ones with just an increase of a constant factor in the number of data samples.  This observation suggests a potentially deeper result along the following lines---many computationally intractable convex sets for which there exist no tight efficiently-computable approximations as measured by integrality gap can nonetheless be well-approximated by computationally tractable convex sets, if the quality of approximation is measured based on statistical inference objectives.

\section*{Acknowledgments}
This material is based upon work supported in part by the U. S. Army Research Laboratory and the U. S. Army Research Office under contract/grant number W911NF-11-1-0391.  We are grateful to Pablo Parrilo, Benjamin Recht, and Parikshit Shah for many insightful conversations.  We would also like to thank Alekh Agarwal, Emmanuel Cand\`{e}s, James Saunderson, Leonard Schulman, and Martin Wainwright for helpful questions and discussions.

\bibliography{refs2}

\newpage

\section*{Supplementary Information}

\subsection*{Proof of Proposition~\ref{prop:denoise2}}

As with the proof of Proposition~\ref{prop:denoise}, we condition on $\bz = \tz$.  Setting $\bdelta = \bx - \tilde{\bx}$ and setting $\hat{\bdelta}_n(\bc) = \hat{\bx}_n(\bc)|_{\bz = \tz} - \tilde{\bx}$, we can rewrite the problem \eqref{eq:shrink} as follows:
\begin{equation*}
\hat{\bdelta}_n(\bc) = \arg \min_{\bdelta \in \R^p} ~~~ \frac{1}{2}\left\|(\bxa - \tilde{\bx}) + \tfrac{\sigma}{\sqrt{n}} \tz - \bdelta \right\|_{\ell_2}^2 ~~~ \mathrm{s.t.} ~~~ \bdelta \in \bc - \tilde{\bx}.
\end{equation*}
Letting $R_1$ and $R_2$ denote orthogonal subspaces that contain $Q_1$ and $Q_2$, i.e., $Q_1 \subseteq R_1$ and $Q_2 \subseteq R_2$, and letting $\bdelta^{(1)} = \mathcal{P}_{R_1}(\bdelta), \bdelta^{(2)} = \mathcal{P}_{R_2}(\bdelta), \hat{\bdelta}^{(1)}_n(\bc) = \mathcal{P}_{R_1}(\hat{\bdelta}_n(\bc)), \hat{\bdelta}^{(2)}_n(\bc) = \mathcal{P}_{R_2}(\hat{\bdelta}_n(\bc))$ denote the projections of $\bdelta, \hat{\bdelta}_n(\bc)$ onto $R_1, R_2$, we can rewrite the above reformulated optimization problem as:
\begin{eqnarray*}
\left[\hat{\bdelta}^{(1)}_n(\bc), \hat{\bdelta}^{(2)}_n(\bc)\right] =  \arg \min_{\bdelta^{(1)} \in Q_1,\bdelta^{(2)} \in Q_2} && ~~ \frac{1}{2}\left\|\mathcal{P}_{R_1}\left[(\bxa - \tilde{\bx}) + \tfrac{\sigma}{\sqrt{n}} \tz\right] - \bdelta^{(1)} \right\|_{\ell_2}^2 \\ && ~~~ + \frac{1}{2}\left\|\mathcal{P}_{R_2}\left[(\bxa - \tilde{\bx}) + \tfrac{\sigma}{\sqrt{n}} \tz\right] - \bdelta^{(2)} \right\|_{\ell_2}^2.
\end{eqnarray*}
As the sets $Q_1,Q_2$ live in orthogonal subspaces, the two variables $\bdelta^{(1)}, \bdelta^{(2)}$ in this problem can be optimized separately.  Consequently, we have that $\|\hat{\bdelta}^{(2)}_n(\bc)\|_{\ell_2} \leq \alpha$ and that
\begin{equation*}
\|\hat{\bdelta}^{(1)}_n(\bc)\|_{\ell_2} \leq \sup_{\bar{\bdelta} \in \mathrm{cone}(Q_1) \cap B^p_{\ell_2}} ~ \langle \bar{\bdelta}, \tfrac{\sigma}{\sqrt{n}} \tz + (\bxa - \tilde{\bx}) \rangle.
\end{equation*}
This bound can be established following the same sequence of steps as in the proof of Proposition~\ref{prop:denoise}.  Combining the two bounds on $\hat{\bdelta}^{(1)}_n(\bc)$ and $\hat{\bdelta}^{(2)}_n(\bc)$, one can then check that
\begin{equation*}
\|\hat{\bdelta}^{(1)}_n(\bc)\|_{\ell_2}^2 + \|\hat{\bdelta}^{(2)}_n(\bc)\|_{\ell_2}^2 \leq 2 \left[ \tfrac{\sigma^2}{n} \G(\mathrm{cone}(Q_1) \cap B^p_{\ell_2}) + \|\bxa - \tilde{\bx}\|_{\ell_2}^2 \right] ~ + ~ \alpha^2.
\end{equation*}
To obtain a bound on $\|\hat{\bx}_n(\bc)|_{\bz = \tz} - \bxa\|_{\ell_2}^2$ we note that
\begin{eqnarray*}
\|\hat{\bx}_n(\bc)|_{\bz = \tz} - \bxa\|_{\ell_2}^2 &\leq& 2 \left[ \|\hat{\bx}_n(\bc)|_{\bz = \tz} - \tilde{\bx}\|_{\ell_2}^2 + \|\bxa - \tilde{\bx}\|_{\ell_2}^2 \right] \\ &\leq& 2\|\hat{\bdelta}^{(1)}_n(\bc)\|_{\ell_2}^2 + 2\|\hat{\bdelta}^{(2)}_n(\bc)\|_{\ell_2}^2 + 2 \|\bxa - \tilde{\bx}\|_{\ell_2}^2.
\end{eqnarray*}
Taking expectations concludes the proof. $\square$

\subsection*{Proof of Proposition~\ref{prop:isop}}
The main steps of this proof follow the steps of a similar result in \cite{ChaRPW2012}, with the principal difference being that we wish to bound Gaussian squared-complexity rather than Gaussian complexity.  A central theme in this proof is the appeal to Gaussian isoperimetry.  Let $\S^{p-1}$ denote the sphere in $p$ dimensions.  Then in bounding the expected squared-distance to the dual cone $\bk^\ast$ with $\bk^\ast \cap \S^{p-1}$ having a volume of $\mu$, we need only consider the extremal case of a spherical cap in $\S^{p-1}$ having a volume of $\mu$.  The manner in which this is made precise will become clear in the proof.  Before proceeding with the main proof, we state and derive a result on the solid angle subtended by a spherical cap in $\S^{p-1}$ to which we will need to appeal repeatedly:
\begin{lemma}\label{lemm:angle}
Let $\psi(\mu)$ denote the solid angle subtended by a spherical cap in $\S^{p-1}$ with volume $\mu \in \left(\tfrac{1}{4} \exp\{-p/20\}, \tfrac{1}{4e^2}\right)$.  Then
\begin{equation*}
\psi(\mu) \geq \frac{\pi}{2}\left(1 - \sqrt{\frac{2 \log\left(\tfrac{1}{4\mu}\right)}{p-1}} \right)
\end{equation*}
\end{lemma}
\noindent \textbf{Proof of Lemma~\ref{lemm:angle}}: Consider the following definition of a spherical cap, parametrized by height $h$:
\begin{equation*}
J = \{\ba \in \S^{p-1} ~| ~ \ba_1 \geq h\}.
\end{equation*}
Here $\ba_1$ denotes the first coordinate of $\ba \in \R^p$.  Given a spherical cap of height $h \in [0, 1]$, the solid angle $\psi$ is given by:
\begin{equation}
\psi = \frac{\pi}{2} - \sin^{-1}(h). \label{eq:angleheight}
\end{equation}
We can thus obtain bounds on the solid angle of a spherical cap via bounds on its height.  The following result from \cite{Bri1998} relates the volume of a spherical cap to its height:
\begin{lemma}\cite{Bri1998}\label{lemm:heightvolume}
For $\tfrac{2}{\sqrt{p}} \leq h \leq 1$ the volume $\tilde{\mu}(p,h)$ of a spherical cap of height $h$ in $\S^{p-1}$ is bounded as
\begin{equation*}
\frac{1}{10 h \sqrt{p}}(1-h^2)^{\tfrac{p-1}{2}} \leq \tilde{\mu}(p,h) \leq \frac{1}{2 h \sqrt{p}}(1-h^2)^{\tfrac{p-1}{2}}
\end{equation*}
\end{lemma}

\noindent Continuing with the proof of Lemma~\ref{lemm:angle}, note that for $\tfrac{2}{\sqrt{p}} \leq h \leq 1$
\begin{equation*}
\frac{1}{2 h \sqrt{p}} (1-h^2)^{\tfrac{p-1}{2}} \leq \frac{1}{4} (1-h^2)^{\tfrac{p-1}{2}} \leq \frac{1}{4} \exp\left( -\tfrac{p-1}{2} h^2 \right).
\end{equation*}
Choosing $h = \sqrt{\tfrac{2 \log\left(\tfrac{1}{4\mu}\right)}{p-1}}$ we have $\tfrac{2}{\sqrt{p}} \leq h \leq 1$ based on the assumption $\mu \in \left(\tfrac{1}{4} \exp\{-p/20\}, \tfrac{1}{4e^2} \right)$.  Consequently, we can apply Lemma~\ref{lemm:heightvolume} with this value of $h$ combined with \eqref{eq:angleheight} to conclude that
\begin{equation*}
\tilde{\mu}\left(p, \sqrt{\frac{2 \log\left(\tfrac{1}{4\mu}\right)}{p-1}}\right) \leq \mu.
\end{equation*}
Hence the solid angle $\psi\left(\tilde{\mu}\left(p, \sqrt{\tfrac{2 \log\left(\tfrac{1}{4\mu}\right)}{p-1}}\right)\right)$ is less than the solid angle $\psi(\mu)$.  Consequently, we use \eqref{eq:angleheight} to conclude that
\begin{equation*}
\psi(\mu) \geq \frac{\pi}{2} - \sin^{-1}\left(\sqrt{\frac{2 \log\left(\tfrac{1}{4\mu}\right)}{p-1}}\right).
\end{equation*}
Using the bound $\sin^{-1}(h) \leq \tfrac{\pi}{2} h$, we obtain the desired bound. $\square$

\textbf{Proof of Proposition~\ref{prop:isop}}: We bound the Gaussian squared-complexity of $\bk$ by bounding the expected squared-distance to the polar cone $\bk^\ast$.  Let $\bar{\mu}(U; t)$ for $U \subseteq \S^{p-1}$ and $t > 0$ denote the volume of the set of points in $\S^{p-1}$ that are within a Euclidean distance of at most $t$ from $U$ (recall that the volume of this set is equivalent to the measure of the set with respect to the normalized Haar measure on $\S^{p-1}$).  We have the following sequence of relations by appealing to the independence of the direction $\bg / \|\bg\|_{\ell_2}$ and of the length $\|\bg\|_{\ell_2}$ of a standard normal vector $\bg$:
\begin{eqnarray*}
\mathbb{E}[\mathrm{dist}(\bg,\bk^\ast)^2] &=& \mathbb{E}[\|\bg\|_{\ell_2}^2 \mathrm{dist}(\bg / \|\bg\|_{\ell_2},\bk^\ast)^2] \\ &=& p ~ \mathbb{E}[\mathrm{dist}(\bg / \|\bg\|_{\ell_2},\bk^\ast)^2] \\ &\leq& p ~ \mathbb{E}[\mathrm{dist}(\bg / \|\bg\|_{\ell_2},\bk^\ast \cap \S^{p-1})^2] \\ &=& p \int_0^\infty \mathbb{P}[\mathrm{dist}(\bg / \|\bg\|_{\ell_2}, \bk^\ast \cap \S^{p-1})^2 > t] dt \\ &=& p \int_0^\infty \mathbb{P}[\mathrm{dist}(\bg / \|\bg\|_{\ell_2}, \bk^\ast \cap \S^{p-1}) > \sqrt{t}] dt \\ &=& 2p \int_0^\infty s \mathbb{P}[\mathrm{dist}(\bg / \|\bg\|_{\ell_2}, \bk^\ast \cap \S^{p-1}) > s] ds \\ &=& 2p \int_0^\infty s [1-\bar{\mu}(\bk^\ast \cap \S^{p-1}; s)] ds.
\end{eqnarray*}
Here the third equality follows based on the integral version of the expected value.  Let $V \subseteq \S^{p-1}$ denote a spherical cap with the same volume $\mu$ as $\bk^\ast \cap \S^{p-1}$.  Then we have by spherical isoperimetry that $\bar{\mu}(V; s) \geq \bar{\mu}(\bk^\ast \cap \S^{p-1}; s)$ for all $s \geq 0$ \cite{Led2000}.  Thus
\begin{equation}
\mathbb{E}[\mathrm{dist}(\bg,\bk^\ast)^2] \leq 2p \int_0^\infty s [1-\bar{\mu}(V; s)] ds. \label{eq:gsc2int}
\end{equation}
From here onward, we focus exclusively on bounding the integral.

Let $\tau(\psi)$ denote the volume of a spherical cap subtending a solid angle of $\psi$ radians.  Recall that $\psi$ is a quantity between $0$ and $\pi$.  As in Lemma~\ref{lemm:angle} let $\psi(\mu)$ denote the solid angle of a spherical cone subtending a solid angle of $\mu$.  Since the Euclidean distance between points on a sphere is always smaller than the geodesic distance, we have that $\bar{\mu}(V; s) \geq \tau(\psi(\mu)+s)$.  Further, we have the following explicit formula for $\tau(\psi)$ \cite{KlaR1997}:
\begin{equation*}
\tau(\psi) = \omega_p^{-1} \int_0^\psi \sin^{p-1}(v)dv,
\end{equation*}
where $\omega_p = \int_0^\pi \sin^{p-1}(v)dv$ is the normalization constant.  Combining these latter two observations, we can bound the integral in \eqref{eq:gsc2int} as:
\begin{eqnarray*}
\int_0^\infty s [1-\bar{\mu}(V; s)] ds &\leq& \int_0^\infty s [1-\tau(\psi(\mu)+s)] ds \\ &=& \int_0^{\pi - \psi(\mu)} s [1-\tau(\psi(\mu)+s)] ds \\ &=& \frac{(\pi-\psi(\mu))^2}{2} - \int_0^{\pi-\psi(\mu)} s \tau(\psi(\mu)+s) ds \\ &=& \frac{(\pi-\psi(\mu))^2}{2} - \omega_p^{-1} \int_0^{\pi-\psi(\mu)} \int_0^{\psi(\mu)+s} s \sin^{p-1}(v)dv ds
\end{eqnarray*}
Next we change the order of integration to obtain:
\begin{eqnarray*}
\int_0^\infty s [1-\bar{\mu}(V; s)] ds &\leq& \frac{(\pi-\psi(\mu))^2}{2} - \omega_p^{-1} \int_0^\pi \int_{\max\{v-\psi(\mu),0\}}^{\pi-\psi(\mu)} \sin^{p-1}(v) s ds dv \\ &=& \frac{(\pi-\psi(\mu))^2}{2} - \omega_p^{-1} \int_0^\pi \frac{1}{2}\left[(\pi-\psi(\mu))^2 - (\max\{v-\psi(\mu),0\})^2 \right] \sin^{p-1}(v) dv \\ &=& \frac{\omega_p^{-1}}{2} \int_0^\pi (\max\{v-\psi(\mu),0\})^2 \sin^{p-1}(v) dv \\ &=& \frac{\omega_p^{-1}}{2} \int_{\psi(\mu)}^\pi (v-\psi(\mu))^2 \sin^{p-1}(v) dv.
\end{eqnarray*}
We now appeal to the inequalities $\omega_p^{-1} \leq \sqrt{p-1} / 2$ and $\sin(x) \leq \exp(-(x-\tfrac{\pi}{2})^2/2)$ for $x \in [0,\pi]$ to obtain
\begin{eqnarray*}
\int_0^\infty s [1-\bar{\mu}(V; s)] ds &\leq& \frac{\sqrt{p-1}}{2} \int_{\psi(\mu)}^\pi (v-\psi(\mu))^2 \exp\left[-\tfrac{p-1}{2}(v - \tfrac{\pi}{2})^2 \right] dv.
\end{eqnarray*}
Performing a change of variables with $a = \sqrt{p-1} (v - \tfrac{\pi}{2})$, we have
\begin{eqnarray*}
\int_0^\infty s [1-\bar{\mu}(V; s)] ds &\leq& \frac{1}{2} \int_{\sqrt{p-1}(\psi(\mu)-\pi/2)}^{\sqrt{p-1} \pi/2} (\tfrac{a}{\sqrt{p-1}} + (\tfrac{\pi}{2} - \psi(\mu)))^2 \exp[-\tfrac{a^2}{2} ] da \\ &\hspace{-2.8in}=& \hspace{-1.4in} \frac{1}{2} \int_{\sqrt{p-1}(\psi(\mu)-\pi/2)}^{\sqrt{p-1} \pi/2} \left[\tfrac{a^2}{p-1} + (\tfrac{\pi}{2} - \psi(\mu))^2 + \tfrac{2a}{\sqrt{p-1}}(\tfrac{\pi}{2} - \psi(\mu)) \right] \exp[-\tfrac{a^2}{2} ] da \\ &\hspace{-2.8in}\leq& \hspace{-1.4in} \frac{1}{2} \left[\int_{-\infty}^{\infty} \tfrac{a^2}{p-1}\exp[-\tfrac{a^2}{2} ] da + \int_{-\infty}^{\infty} (\tfrac{\pi}{2} - \psi(\mu))^2 \exp[-\tfrac{a^2}{2} ] da + \int_{0}^{\infty} \tfrac{2a}{\sqrt{p-1}}(\tfrac{\pi}{2} - \psi(\mu)) \exp[-\tfrac{a^2}{2} ] da \right] \\ &\hspace{-2.8in}=& \hspace{-1.4in} \frac{1}{2} \left[ \tfrac{\sqrt{2\pi}}{p-1} + \sqrt{2\pi}(\tfrac{\pi}{2} - \psi(\mu))^2 + \tfrac{2}{\sqrt{p-1}}(\tfrac{\pi}{2} - \psi(\mu))\cdot (-\exp[-\tfrac{a^2}{2}])|_{0}^{\infty} \right] \\ &\hspace{-2.8in}=& \hspace{-1.4in} \frac{1}{2} \left[  \tfrac{\sqrt{2\pi}}{p-1} + \sqrt{2\pi}(\tfrac{\pi}{2} - \psi(\mu))^2 + \tfrac{2}{\sqrt{p-1}}(\tfrac{\pi}{2} - \psi(\mu))  \right]
\end{eqnarray*}
Here the inequality was obtained by suitably changing the limits of integration.  We now employ Lemma~\ref{lemm:angle} to obtain the final bound:
\begin{eqnarray*}
\G(\bk \cap B_{\ell_2}^p) &\leq& p \left[ \tfrac{\sqrt{2\pi}}{p-1} + \sqrt{2\pi}\left(\tfrac{\pi}{2}\sqrt{\tfrac{2 \log\left(\tfrac{1}{4\mu}\right)}{p-1}}\right)^2 + \tfrac{2}{\sqrt{p-1}} \left(\tfrac{\pi}{2}\sqrt{\tfrac{2 \log\left(\tfrac{1}{4\mu}\right)}{p-1}}\right) \right] \\ &=& \tfrac{p \sqrt{2\pi}}{p-1}  \left[ 1 + \pi \log\left(\tfrac{1}{4 \mu}\right) + \sqrt{\pi} \sqrt{\log\left(\tfrac{1}{4\mu}\right)} \right] \\ &\leq& 20 \log\left(\tfrac{1}{4 \mu}\right).
\end{eqnarray*}
Here the final bound holds because $\mu < 1 / 4e^2$ and $p \geq 12$. $\square$

\end{document}